\documentclass{conm-p-l}

\copyrightinfo{}{}

\newtheorem{theorem}{Theorem}[section]

\theoremstyle{definition}
\newtheorem{definition}[theorem]{Definition}

\theoremstyle{remark}

\newcommand{\GF}[2][2]{{\mathbb F}_{#1^{#2}}}
\newcommand{\tr}[2][1]{{\rm Tr}_{#1}^{#2}}

\numberwithin{equation}{section}

\begin{document}

\title{Niho Bent Functions and Subiaco/Adelaide Hyperovals}


\author{Tor Helleseth}
\address{Department of Informatics, University of Bergen, P.O. Box 7800, N-5020 Bergen, Norway}
\email{Tor.Helleseth@ii.uib.no}

\author{Alexander Kholosha}
\address{Department of Informatics, University of Bergen, P.O. Box 7800, N-5020 Bergen, Norway}
\email{Alexander.Kholosha@ii.uib.no}

\author{Sihem Mesnager}
\address{Department of Mathematics, University of Paris 8 and University of Paris 13, 2 rue de la
libert\'e, 93526 Saint-Denis Cedex, France}
\email{smesnager@univ-paris8.fr}

\subjclass[2010]{Primary }

\date{}

\begin{abstract}
In this paper, the relation between binomial Niho bent functions
discovered by Dobbertin {\em et al.} and o-polynomials that give rise
to the Subiaco and Adelaide classes of hyperovals is found. This
allows to expand the class of bent functions that corresponds to
Subiaco hyperovals, in the case when $m\equiv 2\ (\bmod\;4)$.
\end{abstract}

\maketitle

\section{Introduction and Preliminaries}
Boolean functions of $n$ variables are binary functions over the
Galois field $\GF n$ (or over the vector space ${\mathbb F}_2^n$ of
all binary vectors of length $n$). In this paper, we shall always
endow this vector space with the structure of a field, thanks to the
choice of a basis of $\GF n$ over $\GF {}$. Boolean functions are
used in the pseudo-random generators of stream ciphers and play a
central role in their security.

Bent functions were introduced by Rothaus \cite{Ro76} in 1976. These
are Boolean functions of an even number of variables $n$, that are
maximally nonlinear in the sense that their Walsh transform takes
precisely the values $\pm 2^{n/2}$. This corresponds to the fact that
their Hamming distance to all affine functions is optimal. Bent
functions have also attracted a lot of research interest because of
their relations to coding theory and applications in cryptography.
Despite their simple and natural definition, bent functions turned
out to admit a very complicated structure in general. On the other
hand, many special explicit constructions are known. Distinguished
are primary constructions giving bent functions from scratch and
secondary ones building new bent functions from one or several given
bent functions. These constructions often look simpler when written
in their bivariate representation but, of course, they also have an
equivalent univariate form (see Subsection~\ref{subsec:unibi}).

It is well known that some of the explicit constructions belong to
the two general families of bent functions which are the original
Maiorana-McFarland \cite{McFa73} and the Partial Spreads classes. It
was in the early seventies when Dillon in his thesis \cite{Di74}
introduced the two above mentioned classes plus the third one denoted
by $H$, where bentness is proven under some conditions which were not
obvious to achieve (in this class, Dillon was able to exhibit only
those functions belonging, up to the affine equivalence, to the
Maiorana-McFarland class). He defined the functions in class $H$ in
their bivariate representation but they can also be seen in the
univariate form as found recently by Carlet and Mesnager
\cite{CaMe11}. They extended the class $H$ to a slightly larger class
denoted by $\mathcal{H}$ (see Subsection~\ref{subsec:class_H}), also
defined in bivariate representation, and observed that this class
contains all bent functions of the, so called, Niho type which were
introduced in \cite{DoLeCaCaFeGa06} by Dobbertin {\em et al.} (see
Subsection~\ref{subsec:Niho}).

\subsection{Trace representation, Boolean functions in univariate
and bivariate forms.}
 \label{subsec:unibi}
For any positive integer $k$ and any $r$ dividing $k$, the trace
function $\tr[r]k()$ is the mapping from $\GF k$ to $\GF r$ defined
by
\[\tr[r]k(x):=\sum_{i=0}^{\frac{k}{r}-1}x^{2^{ir}}=x+x^{2^r}+x^{2^{2r}}+\cdots+x^{2^{k-r}}\enspace.\]
In particular, the {\em absolute trace} over $\GF k$ is the function
$\tr k(x)=\sum_{i=0}^{k-1}x^{2^i}$. Recall that the trace function
satisfies the transitivity property $\tr k=\tr r\circ\tr[r]k$. From
now on assume $n$ is even and $n=2m$. For any set $E$, denote
$E\setminus\{0\}$ by $E^*$.

The {\em univariate representation} of a Boolean function is defined
as follows: we identify ${\mathbb F}_2^n$ (the $n$-dimensional vector
space over $\GF {}$) with $\GF n$ and consider the arguments of $f$
as elements in $\GF n$. An inner product in $\GF n$ is $x\cdot y=\tr
n(xy)$. There exists a unique univariate polynomial
$\sum_{i=0}^{2^n-1}a_i x^i$ over $\GF n$ that represents $f$ (this is
true for any vectorial function from $\GF n$ to itself). The
algebraic degree of $f$ is equal to the maximum $2$-weight of an
exponent having nonzero coefficient, where the $2$-weight $w_2(i)$ of
an integer $i$ is the number of ones in its binary expansion. Hence,
in the case of a bent function, all exponents $i$ whose $2$-weight is
larger than $m$ have a zero coefficient $a_i$. Moreover, $f$ being
Boolean, its univariate representation can be written in the form of
$f(x)=\sum_{j\in\Gamma_n}\tr {o(j)}(a_j x^j)$, where $\Gamma_n$ is
the set of integers obtained by choosing one element in each
cyclotomic coset of $2$ modulo $2^n-1$, $o(j)$ is the size of the
cyclotomic coset containing $j$ and $a_j\in\GF {o(j)}$. This
representation is unique up to the choice of cyclotomic coset
representatives. Function $f$ can also be written in a non-unique way
as $\tr n(P(x))$ where $P(x)$ is a polynomial over $\GF n$.

The {\em bivariate representation} of a Boolean function is defined
as follows: we identify $\GF n$ with $\GF m\times\GF m$ and consider
the argument of $f$ as an ordered pair $(x,y)$ of elements in $\GF
m$. There exists a unique bivariate polynomial $\sum_{0\leq i,j\leq
2^m-1}a_{i,j}x^i y^j$ over $\GF m$ that represents $f$. The algebraic
degree of $f$ is equal to
\[\max_{(i,j)\,|\,a_{i,j}\neq 0}(w_2(i)+w_2(j))\enspace.\]
And $f$ being Boolean, its bivariate representation can be written in
the form of $f(x,y)=\tr m(P(x,y))$, where $P(x,y)$ is some polynomial
of two variables over $\GF m$.

Let $f$ be an $n$-variable Boolean function. Its {\em ``sign"
function} is the integer-valued function $\chi_f:=(-1)^f$. The {\em
Walsh transform} of $f$ is the discrete Fourier transform of $\chi_f$
whose value at point $w\in\GF n$ is defined by
\[\hat\chi_f(w)=\sum_{x\in\GF n}(-1)^{f(x)+\tr n(wx)}\enspace.\]

\begin{definition}
For even $n$, a Boolean function $f$ in $n$ variables is said to be
bent if for any $w\in\GF n$ we have $\hat\chi_f(w)=\pm
2^{\frac{n}{2}}$.
\end{definition}

\subsection{Class $\mathcal{H}$ of Bent Functions}
 \label{subsec:class_H}
In his thesis \cite{Di74}, Dillon introduced the class of bent
functions denoted by $H$. The functions in this class are defined in
their bivariate form as
\[f(x,y)=\tr m\big(y+xG(yx^{2^m-2})\big)\enspace,\]
where $x,y\in\GF m$ and $G$ is a permutation of $\GF m$ such that
$G(x)+x$ does not vanish and for any $\beta\in\GF m^*$, the function
$G(x)+\beta x$ is $2$-to-$1$ (i.e., the pre-image of any element of
$\GF m$ is either a pair or the empty set). As observed by Carlet and
Mesnager \cite[Proposition~1]{CaMe11}, this class can be slightly
extended into a class $\mathcal{H}$ defined as the set of (bent)
functions $g$ satisfying
 \setlength{\arraycolsep}{0.14em}
\begin{equation}
 \label{eq:H_class}
g(x,y)=\left\{\begin{array}{ll}
\tr m\left(xH\left(\frac{y}{x}\right)\right),&\ \mbox{if}\quad x\neq 0\\
\tr m(\mu y),&\ \mbox{if}\quad x=0\enspace,
\end{array}\right.
\end{equation}
 \setlength{\arraycolsep}{5pt}\noindent
where $\mu\in\GF m$ and $H$ is a mapping from $\GF m$ to itself
satisfying the following necessary and sufficient conditions
\begin{align}
 \label{eq:co1}
G:z&\mapsto H(z)+\mu z\ \mbox{is a permutation on}\ \GF m\\
 \label{eq:co2}
z&\mapsto G(z)+\beta z\ \mbox{is 2-to-1 on}\ \GF m\ \mbox{for any}\ \beta\in\GF m^*\enspace.
\end{align}
As proved in \cite[Lemma~13]{CaMe11}, condition (\ref{eq:co2})
implies condition (\ref{eq:co1}) and, thus, is necessary and
sufficient for $g$ being bent. It also follows that polynomials
$G(z)$ satisfying (\ref{eq:co2}) are so-called o-polynomials (oval
polynomials) over $\GF m$ (the additional properties of $G(0)=0$ and
$G(1)=1$ can be achieved by taking $\frac{G(z)+G(0)}{G(1)+G(0)}$
instead of $G(z)$). o-polynomials arise from hyperovals and define
them. Note that class ${\mathcal H}$ contains all bent functions with
the property that their restriction to the multiplicative cosets of
$\GF m$ is linear.

\subsection{Niho bent functions}
 \label{subsec:Niho}
Recall that a positive integer $d$ (always understood modulo $2^n-1$)
is said to be a {\em Niho exponent} and $t\mapsto t^d$ is a {\em Niho
power function} if the restriction of $t^d$ to $\GF m$ (and,
therefore, to its cosets $u\GF m$) is linear or, in other words,
$d\equiv 2^j\ (\bmod\;2^m-1)$ for some $j<n$. As we consider $\tr[1]
n(at^d)$ with $a\in\GF n$, without loss of generality, we can assume
that $d$ is in the normalized form, i.e., with $j=0$. Then we have a
unique representation $d=(2^m-1)s+1$ with $2\leq s\leq 2^m$. The
simplest example of an infinite class of Niho bent functions is the
quadratic function $\tr[1] m(at^{2^m +1})$ with $a\in\GF m^*$. Other
known classes are:
\begin{itemize}
\item Three examples from \cite{DoLeCaCaFeGa06} are binomials of
    the form $f(t)=\tr[1] n(\alpha_1 t^{d_1}+ \alpha_2 t^{d_2})$,
    where $2d_1=2^m+1\in{\mathbb Z}/(2^n-1){\mathbb Z}$ and
    $\alpha_1,\alpha_2\in\GF n^*$ are such that
    $(\alpha_1+\alpha_1^{2^m})^2=\alpha_2^{2^m+1}$. Equivalently,
    denoting $a=(\alpha_1+\alpha_1^{2^m})^2$ and $b=\alpha_2$ we
    have $a=b^{2^m+1}\in\GF m^*$ and $f(t)=\tr[1]
    m(at^{2^m+1})+\tr[1] n(bt^{d_2})$. Note that if $b=0$ and
    $a\neq 0$ then $f$ is also bent but becomes quadratic equal
    to the function mentioned above. The possible values of $d_2$
    are:
\begin{itemize}
\item[] $d_2=(2^m-1)3+1$ (with the condition that, if
    $m\equiv 2\ (\bmod\;4)$ then $b$ is the fifth power of an
    element in $\GF n$; otherwise, $b$ can be any nonzero
    element),
\item[] $4d_2=(2^m-1)+4$ (with the condition that $m$ is
    odd),
\item[] $6d_2=(2^m-1)+6$ (with the condition that $m$ is
    even).
\end{itemize}
As observed in \cite{DoLeCaCaFeGa06,CaHeKhMe11}, these functions
have algebraic degree $m$, 3 and $m$ respectively.
\item An extension by Leander and Kholosha \cite{LeKh06} of the
    second class from \cite{DoLeCaCaFeGa06} has the form of
\begin{equation}
 \label{eq:LKh_bent}
\tr[1] n\Big(at^{2^m+1}+\sum_{i=1}^{2^{r-1}-1}t^{(2^m-1)\frac{i}{2^r}+1}\Big)
\end{equation}
    with $r>1$ satisfying $\gcd(r,m)=1$ and $a\in\GF n$ is such
    that $a+a^{2^m}=1$.
\item Functions in a bivariate form obtained from the known
    o-polynomials (see \cite[Section~6]{CaMe11}).
\end{itemize}
As was noted in \cite{DoLeCaCaFeGa06}, all cases except for
$d_2=(2^m-1)3+1$ with $m\equiv 2\ (\bmod\;4)$ give
$\gcd(d_2,2^n-1)=1$ and in the remaining case, $\gcd(d_2,2^n-1)=5$.
Therefore, having the condition on $b$, it can be assumed, without
loss of generality, that $b=1$ (this is achieved by substituting $t$
with $b^{-1/d_2}t$). However, in Subsection~\ref{subsec:2}, we show
that even in the case when $m\equiv 2\ (\bmod\;4)$ the value of $b$
can be taken arbitrary under the condition that $a=b^{2^m+1}$.

Since the restriction to $u\GF m$ of these bent functions is linear,
they all belong to ${\mathcal H}$. The question left open in
\cite{DoLeCaCaFeGa06} was finding the dual and checking if that was
of the Niho type (possibly up to affine equivalence). In
\cite{CaMe11,CaHeKhMe11} considered were bent functions from the
second class (having degree $3$) and multinomial (\ref{eq:LKh_bent}).
It was shown that corresponding o-polynomials are Frobenius mappings
and dual functions were found that turned out not to be in the Niho
class. Moreover, these cases give bent functions in the completed
Maiorana-McFarland class. In this paper, we find o-polynomials that
arise from the first class of binomial Niho bent functions. However,
it still remains to determine the dual.

\section{Subiaco Hyperovals}
Here we define o-polynomials that give rise to the Subiaco family of
hyperovals.

\begin{theorem}[Theorems~3-5~\cite{ChPePiRo96}]
 \label{th:Sub}
Take polynomials $f(x)$ and $g(x)$ and for any $s\in\GF m$ define
\begin{equation}
 \label{eq:fs}
f_s(x)=\frac{f(x)+esg(x)+s^{1/2}x^{1/2}}{1+es+s^{1/2}}\enspace,
\end{equation}
where $e\in\GF m$ with $\tr m(e)=1$ is defined further. Then in the
following cases, $g(x)$ and $f_s(x)$ are o-polynomials:
\renewcommand{\theenumi}{\roman{enumi}}
\renewcommand{\labelenumi}{(\theenumi)}
\begin{enumerate}
\item\label{it:1} if $m$ is odd then take $e=1$ and
\[f(x)=\frac{x^2+x}{(x^2+x+1)^2}+x^\frac{1}{2}\quad\mbox{and}\quad g(x)=\frac{x^4+x^3}{(x^2+x+1)^2}+x^\frac{1}{2}\enspace;\]
\item\label{it:2} if $m\equiv 2\ (\bmod\;4)$ then take $e=w\in\GF
    m$ with $w^2+w+1=0$ and
\[f(x)=\frac{x^2(x^2+wx+w)}{(x^2+wx+1)^2}+w^2 x^\frac{1}{2}\quad\mbox{and}\quad
  g(x)=\frac{wx(x^2+x+w^2)}{(x^2+wx+1)^2}+w^2 x^\frac{1}{2}\enspace;\]
\item\label{it:3} for any $m$, take
    $e=\frac{w^2+w^5+w^{1/2}}{w(1+w+w^2)}$ where $w\in\GF m$ with
    $w^2+w+1\neq 0$ and $\tr m(1/w)=1$, and
\[\begin{split}
f(x)&=\frac{w^2(x^4+x)+w^2(1+w+w^2)(x^3+x^2)}{(x^2+wx+1)^2}+x^\frac{1}{2}\quad\mbox{and}\\
g(x)&=\frac{w^4 x^4+w^3(1+w^2+w^4)x^3+w^3(1+w^2)x}{(w^2+w^5+w^{1/2})(x^2+wx+1)^2}+\frac{w^{1/2}}{w^2+w^5+w^{1/2}}x^\frac{1}{2}\enspace.
\end{split}\]
\end{enumerate}
\end{theorem}

It is useful to have the following explicit expressions for $f_s(x)$
in each of the cases considered. Denote $1+es+s^\frac{1}{2}=A$, then
$f_s(x)$ is equal to
\begin{align}
\label{eq:m_odd}&\frac{s(x^4+x^3)+x^2+x}{A(x^2+x+1)^2}+x^\frac{1}{2}\ ,\quad m\ \mbox{odd}\\
\label{eq:m/2_odd}&A^{-1}\left(\frac{x^4+w(sw+1)(x^3+x^2)+swx}{(x^2+wx+1)^2}+(w^2+s+s^\frac{1}{2})x^\frac{1}{2}\right)\ ,\quad m/2\ \mbox{odd}\\
\label{eq:m_arb}&\bigg(w^2\frac{(1+sw+w^2)x^4+(1+w+w^2)^2(sx^3+x^2)+(s+w+s w^2)x}{(1+w+w^2)(x^2+wx+1)^2}\\
\nonumber&\quad+\left(s^\frac{1}{2}+\frac{s+1}{w^{1/2}(1+w+w^2)}\right)x^\frac{1}{2}\bigg)(e+es+s^\frac{1}{2})^{-1}\ ,\quad m\ \mbox{arbitrary}\enspace,
\end{align}
where in (\ref{eq:m_arb}), we changed $s+1$ for $s$ in the original
definition of $f_s(x)$. Note that for $m$ odd, taking $w=1$ in
(\ref{eq:m_arb}) results in (\ref{eq:m_odd}).

In each of the cases listed above, the set $(f(x),g(x),a)$ defines a
$q$-clan. On the other hand, by \cite[Theorem~1]{ChPePiRo96}, the
existence of the $q$-clan is equivalent to the property that $g(x)$
is an o-polynomial and $f_s(x)$ is an o-polynomial for any $s\in\GF
m$. In \cite{PaPePi95}, it was shown that the Subiaco construction
provides two inequivalent hyperovals if $m\equiv 2\ (\bmod\;4)$ and
one hyperoval otherwise.

\section{Bent Functions from Subiaco Hyperovals}
Take the following function over $\GF n$
\[f(t)=\tr m(a t^{2^m+1})+\tr n(b t^{3(2^m-1)+1})\enspace,\]
where $a\in\GF m^*$ and $b\in\GF n^*$ are such that $b^{2^m+1}=a$.
Let $(u,v)$ be a basis of $\GF n$ as a two-dimensional vector space
over $\GF m$. Then for any $x,y\in\GF m$, we obtain $f(ux+vy)$ having
the form of (\ref{eq:H_class}) with
\[\begin{split}
H(z)&=a^{\frac{1}{2}}(u+vz)^{\frac{2^m+1}{2}}+\tr[m]n\big(b(u+vz)^{3(2^m-1)+1}\big)\\
\mu&=a^{\frac{1}{2}}v^{\frac{2^m+1}{2}}+\tr[m]n(bv^{3(2^m-1)+1})\enspace.
\end{split}\]
Here all notation are from Subsection~\ref{subsec:class_H}.
Therefore, with $z\in\GF m$,
\[G(z)=a^{\frac{1}{2}}v^{\frac{2^m+1}{2}}z+a^{\frac{1}2}(u+vz)^{\frac{2^m+1}{2}}+
\tr[m]n\big(b(v^{3(2^m-1)+1}z+(u+vz)^{3(2^m-1)+1})\big)\enspace.\]

Further, we have that
\[(u+vz)^{\frac{2^m+1}{2}}=u^{\frac{2^m+1}{2}}+
\big(\tr[m]n(u^{2^{m}}v)\big)^{\frac{1}{2}}z^{\frac{1}{2}}+(vz)^{\frac{2^m+1}{2}}\]
and since $z\in\GF m$,
\begin{equation}
 \label{eq:aux}
a^{\frac{1}{2}}v^{\frac{2^m+1}{2}}z+a^{\frac{1}{2}}(u+vz)^{\frac{2^m+1}{2}}=
a^{\frac{1}{2}}u^{\frac{2^m+1}{2}}+a^{\frac{1}{2}}\big(\tr[m]n(u^{2^{m}}v)\big)^{\frac{1}{2}}z^{\frac{1}{2}}\enspace.
\end{equation}

Now expand the term $(u+vz)^{3(2^m-1)+1}$. To this end, note that
$3(2^m-1)+1=2^{m+1}-1+2^m-1$. Then
\[\begin{split}
(u+vz)^{3(2^m-1)+1}&=(u+vz)^{2^{m+1}-1}(u+vz)^{2^m-1}\\
&=\sum_{j=0}^{2^{m+1}-1}u^{2^{m+1}-1-j}(vz)^j\sum_{j=0}^{2^{m}-1}u^{2^m-1-j}(vz)^j\\
&=\sum_{i=0}^{3\cdot 2^m-2}(N_i\bmod 2)u^{3\cdot2^m-2-i}(vz)^i\enspace,
\end{split}\]
where $N_i=|E_i|$ and
\[E_i=\{(j_1,j_2)\mid j_1+j_2=i,\ 0\leq j_1\leq 2^{m+1}-1,\
0\leq j_2\leq 2^m-1\}\enspace.\] We compute $N_i$ by enumerating the
elements of $E_i$ as follows:
\begin{itemize}
\item for $0\leq i\leq 2^m-1$, we have $E_i=\{(i-j,j)\mid 0\leq
    j\leq i\}$ and $N_i=i+1$;
\item for $2^m\leq i\leq 2^{m+1}-1$, we have $E_i=\{(i-j,j)\mid
    0\leq j\leq 2^m-1\}$ and $N_i=2^m$;
\item for $2^{m+1}\leq i\leq 3\cdot 2^m-2$, we have
    $E_i=\{(i-j,j)\mid i-2^{m+1}+1\leq j\leq 2^m-1\}$ and
    $N_i=3\cdot 2^m-1-i$ (indeed, $j_1+j_2=i$ implies that
    $j_2=i-j_1\geq i-2^{m+1}+1$ since $j_1\leq 2^{m+1}-1$).
\end{itemize}

Therefore $N_i\bmod 2=1$ if and only if $i=2l$ with $0\leq l\leq
2^{m-1}-1$ or $i=2^{m+1}+2l$ with $0\leq l\leq 2^{m-1}-1$ and
\[\begin{split}
(u&+vz)^{3(2^m-1)+1}=\sum_{l=0}^{2^{m-1}-1}u^{3\cdot 2^m-2-2l}(vz)^{2l}+
\sum_{l=0}^{2^{m-1}-1}u^{3\cdot 2^m-2-2^{m+1}-2l}(vz)^{2^{m+1}+2l}\\
&\stackrel{(\ast)}{=}\sum_{l=0}^{2^{m-1}-1}u^{3\cdot 2^m-2(l+1)}(vz)^{2l}+
\sum_{l=0}^{2^{m-1}-1}u^{2^m-2(l+1)}v^{2^{m+1}-2}(vz)^{2(l+1)}\\
&=\sum_{l=0}^{2^{m-1}-1}u^{3\cdot 2^m-2(l+1)}(vz)^{2l}+\sum_{l=1}^{2^{m-1}}u^{2^m-2l}v^{2^{m+1}-2}(vz)^{2l}\\
&=u^{3\cdot 2^m-2}+(u^{3\cdot 2^m-2}+u^{2^m}v^{2^{m+1}-2})\sum_{l=1}^{2^{m-1}-1}\left(u^{-1}vz\right)^{2l}+
v^{3\cdot 2^m-2}z\\
&=u^{3\cdot 2^m-2}+u^{2^m}(u^{2(2^m-1)}+v^{2(2^m-1)})\left(1+\frac{1+(u^{-1}vz)^{2^m}}{1+u^{-2}v^2 z^2}\right)+
v^{3\cdot 2^m-2}z\\
&=u^{2^m}v^{2(2^m-1)}+u^{2^m}(u^{2(2^m-1)}+v^{2(2^m-1)})(1+u^{-1}vz)^{2^m-2}+v^{3\cdot 2^m-2}z\\
&=u^{2^m}v^{2(2^m-1)}+u^2(u^{2(2^m-1)}+v^{2(2^m-1)})(u+vz)^{2^m-2}+v^{3\cdot 2^m-2}z\enspace.
\end{split}\]
In the second sum after ($\ast$), we used that
$z^{2^{m+1}+2l}=(z^{2^{m}})^2 z^{2l}=z^2 z^{2l}=z^{2(l+1)}$. Finally,
denoting
\[c=a^{\frac{1}{2}}u^{\frac{2^m+1}{2}}+\tr[m]n(bu^{2^m}v^{2(2^m-1)})\]
and using (\ref{eq:aux}), we obtain that
\begin{equation}
 \label{eq:G_gen}
G(z)=c+a^{\frac{1}{2}}\big(\tr[m]n(u^{2^{m}}v)\big)^{\frac{1}{2}}z^{\frac{1}{2}}+
\tr[m]n\big(bu^2(u^{2(2^m-1)}+v^{2(2^m-1)})(u+vz)^{2^m-2}\big)\enspace.
\end{equation}

Now assume $v=1$ and take $u\in\GF n\setminus\{1\}$ with
$u^{2^m+1}=1$ that means $u\in\GF n\setminus\GF m$. Also denote
$u+u^{2^m}=w\in\GF m^*$ and observe that $\tr m(1/w)=1$ (since this
is equivalent to $u^2+wu+1$ being irreducible over $\GF m$).
Moreover, all $w\in\GF m$ with such a trace property are obtained in
this way from $u$. Then $u^{2^m-1}=w/u+1$ and
\[\begin{split}
\tr[m]n(u^{2^m} v)&=w\\
u^2\big(v^{2(2^m-1)}+u^{2(2^m-1)}\big)&=w^2\enspace.
\end{split}\]
Under these conditions, $c=a^\frac{1}{2}+\tr[m]n(b u^{2^m})$ and
\begin{align}
 \label{eq:G}
&G(z)=c+(awz)^\frac{1}{2}+\frac{bw^2(u^{2^m}+z)}{(u+z)^2}+\frac{b^{2^m}w^2(u+z)}{(u^{2^m}+z)^2}\\
\nonumber&=c+(awz)^\frac{1}{2}+w^2\frac{b(u+w+z)^3+b^{2^m}(u+z)^3}{(u+z)^2(u+w+z)^2}\\
\nonumber&=c+(awz)^\frac{1}{2}+w^2\frac{(b+b^{2^m})(u+z)^3+bw(z^2+wz+u^{2^m+1}+w^2)}{(z^2+wz+u^{2^m+1})^2}\\
\nonumber&\stackrel{(\ref{eq:1})}{=}c+(awz)^\frac{1}{2}\\
\nonumber&\quad+\frac{w^2(b+b^{2^m})(z^3+uz^2+u^2 z)+bw^3(z^2+wz)+\tr[m]n(b^{2^m}(u^5+u))}{(z^2+wz+1)^2}\\
\nonumber&=a^\frac{1}{2}+\tr[m]n(b^{2^m}u^5)+(awz)^\frac{1}{2}\\
\nonumber&\quad+\frac{w^2(b+b^{2^m})(z^3+uz^2+u^2 z)+bw^3(z^2+wz)+\tr[m]n(b^{2^m}(u^5+u))(z^2+wz)^2}{(z^2+wz+1)^2}\\
\nonumber&\stackrel{(\ref{eq:2},\ref{eq:3})}{=}a^\frac{1}{2}+\tr[m]n(b^{2^m}u^5)+(awz)^\frac{1}{2}\\
\nonumber&\quad+\frac{\tr[m]n(b^{2^m}(u^5+u))z^4+\tr[m]n(b)w^2
z^3+\tr[m]n(b^{2^m}u^5)w^2 z^2+\tr[m]n(b^{2^m}(u^4+1))z}
{(z^2+wz+1)^2}\ .
\end{align}
Here we used the following identities
\begin{align}
\label{eq:1}w^2(b+b^{2^m})u^3+bw^3(1+w^2)&=\tr[m]n(b^{2^m}(u^5+u))\ ;\\
\label{eq:2}u(b+b^{2^m})+bw+\tr[m]n(b^{2^m}(u^5+u))&=\tr[m]n(b^{2^m}u^5)\ ;\\
\label{eq:3}w^2(b+b^{2^m})u^2+bw^4&=\tr[m]n(b^{2^m}(u^4+1))\enspace.
\end{align}

Further, we consider three separate cases defined by the value of
$m$.

\subsection{$m$ odd}
In this case, take $u\in{\mathbb F}_4\setminus\{0,1\}$. Note that
$u\in\GF n\setminus\GF m$ and $w=u+u^{2^m}=u+u^2=1$. Then, by
(\ref{eq:G}),
\[\begin{split}
G(z)&=a^\frac{1}{2}+\tr[m]n(bu)+(az)^\frac{1}{2}+\frac{\tr[m]n(b)(z^4+z^3)+\tr[m]n(bu)(z^2+z)}{(z^2+z+1)^2}\\
&=a^\frac{1}{2}+\tr[m]n(bu)+(az)^\frac{1}{2}+a^\frac{1}{2}\frac{(B+B^{-1})(z^4+z^3)+(B^{-1}u^2+Bu)(z^2+z)}{(z^2+z+1)^2}\\
&=a^\frac{1}{2}+\tr[m]n(bu)+a^\frac{1}{2}f_s(z)\enspace,
\end{split}\]
where $B=ba^{-\frac{1}{2}}$ with
$B^{-1}=b^{2^m}a^{-\frac{1}{2}}=B^{2^m}$ since $a=b^{2^m+1}$.
Polynomial $f_s(z)$ with $s=\frac{1+B^2}{u^2+B^2 u}\in\GF m$ is an
o-polynomial (\ref{eq:m_odd}) (assuming $u^2+B^2 u\neq 0$). In the
case when $u^2=B^2 u$ (or, equivalently, $b^{2^m-1}=u^2$) we obtain
\[G(z)=bu+buz^\frac{1}{2}+bu\frac{z^4+z^3}{(z^2+z+1)^2}=bu(1+g(z))\enspace,\]
since $a^\frac{1}{2}=(b^{2^m+1})^\frac{1}{2}=bu=b+b^{2^m}$ and where
o-polynomial $g(z)$ comes from
Theorem~\ref{th:Sub}~Item~(\ref{it:1}).

Assuming $b^{2^m-1}\neq u^2$, note that equation
$s=\frac{b^{2^m-1}+1}{b^{2^m-1}u^2+u}$ can be solved for the unknown
$b\in\GF n^*$ for any $s\in\GF m$ since $s\neq u$. We conclude that
the set of bent functions with $b\in\GF n^*$ corresponds exactly to
all o-polynomials described in
Theorem~\ref{th:Sub}~Item~(\ref{it:1}). This means that the existence
of this set of bent functions is equivalent to the existence of the
corresponding $q$-clan.

\subsection{$m\equiv 2\ (\bmod\;4)$}
 \label{subsec:2}
In this case, take $u\in{\mathbb F}_{16}\setminus{\mathbb F}_4$ with
$u^5=1$. Note that $u\in\GF n\setminus\GF m$ and $u^{2^m+1}=u^5=1$.
Then $u+u^{2^m}=u+u^4=w\in{\mathbb F}_4\subset\GF m$. Obviously,
$w\neq 0$. It can be checked directly that $u$ with the prescribed
properties also satisfies $w\neq 1$ and, thus, $w^2+w=1$. There are
{\em four} options for choosing $u$ with these properties and both
$w\in{\mathbb F}_4\setminus\{0,1\}$ can be obtained. Then, by
(\ref{eq:G}),
\begin{align*}
&G(z)=a^\frac{1}{2}+\tr[m]n(b)+(awz)^\frac{1}{2}\\
&\quad+\frac{\tr[m]n(b(u^4+1))z^4+\tr[m]n(b)w^2(z^3+z^2)+\tr[m]n(b(u+1))z}{(z^2+wz+1)^2}\\
&=a^\frac{1}{2}+\tr[m]n(b)+(awz)^\frac{1}{2}+\tr[m]n(b(u^4+1))\frac{z^4+w(sw+1)(z^3+z^2)+swz}{(z^2+wz+1)^2}\\
&\stackrel{(\ast)}{=}a^\frac{1}{2}+\tr[m]n(b)+(1+ws+s^\frac{1}{2})\tr[m]n(b(u^4+1))f_s(z)\enspace,
\end{align*}
where polynomial $f_s(z)$ with
$s=\frac{w^2\tr[m]n(b(u+1))}{\tr[m]n(b(u^4+1))}$ is an o-polynomial
(\ref{eq:m/2_odd}) (assuming $\tr[m]n(b(u^4+1))\neq 0$). In the case
when $\tr[m]n(b(u^4+1))=0$ (or, equivalently,
$b^{2^m-1}=(u+1)^3=u^4$) we obtain
\[\begin{split}
G(z)&=a^\frac{1}{2}+\tr[m]n(b)+(awz)^\frac{1}{2}+\frac{\tr[m]n(b)w^2(z^3+z^2)+\tr[m]n(b(u+1))z}{(z^2+wz+1)^2}\\
&=a^\frac{1}{2}+\tr[m]n(b)+b u^2 w^2 z^\frac{1}{2}+b u^2\frac{wz(z^2+z+w^2)}{(z^2+wz+1)^2}\\
&=a^\frac{1}{2}+\tr[m]n(b)+b u^2 g(z)\enspace,
\end{split}\]
since $a=b^{2^m+1}=b^2 u^4$ and $\tr[m]n(b)w=b(1+u^4)(u+u^4)=bu^2$
and where o-polynomial $g(z)$ comes from
Theorem~\ref{th:Sub}~Item~(\ref{it:2}). On the other hand, if
$b^{2^m-1}=u^4$ then it suffices just to take another $u$ with the
above defined properties (recall that four options exist). To obtain
$(\ast)$ we used the following identities
\begin{align*}
&(w+s^2+s)\tr[m]n(b(u^4+1))^2\\
&=w\tr[m]n(b(u^4+1))^2+w\tr[m]n(b(u+1))^2+w^2\tr[m]n(b(u+1))\tr[m]n(b(u^4+1))\\
&=w^2\big(\tr[m]n(bu)\tr[m]n(bu^4)+\tr[m]n(b)\tr[m]n(b(u^4+u))+\tr[m]n(b)^2\big)+w\tr[m]n(b(u^4+u))^2\\
&=w^2(bu+b^{2^m}u^4)(bu^4+b^{2^m}u)+w^2\tr[m]n(b)^2=aw\enspace.
\end{align*}

It is important to observe that there are no restrictions on the
value of $b$ here. It means that this technique allows to enlarge the
original class of Niho bent functions proved in
\cite{DoLeCaCaFeGa06}.

Assuming $b^{2^m-1}\neq u^4$, note that equation
$s=\frac{w^2\tr[m]n(b(u+1))}{\tr[m]n(b(u^4+1))}$ can be solved for
the unknown $b\in\GF n^*$ for any $s\in\GF m$. Indeed, this equation
can be rewritten as
\[\begin{split}
b(u^4 s+s+u w^2+w^2)&=b^{2^m}(us+s+u^4 w^2+w^2)\quad\mbox{or}\\
b(u^4 s+s+u^4+u^2)&=b^{2^m}(us+s+u^3+u)\enspace.
\end{split}\]
Since $s\in\GF m$, it is easy to see that this equation has nonzero
sides and its right-hand side is a $2^m$th power of the left-hand
side. We conclude that the set of bent functions with $b\in\GF n^*$
corresponds exactly to all o-polynomials described in
Theorem~\ref{th:Sub}~Item~(\ref{it:2}). This means that the existence
of this set of bent functions is equivalent to the existence of the
corresponding $q$-clan.

\subsection{$m\equiv 0\ (\bmod\;4)$}
In this case, $w^2+w+1\neq 0$ since the opposite is equivalent to
$u^4+u^3+u^2+u+1=0$ that gives $u\in{\mathbb F}_{2^4}$ which is a
contradiction because ${\mathbb F}_{2^4}\subset\GF m$. As was noted
in Subsection~\ref{subsec:Niho}, without loss of generality, we can
assume $b=a=1$. Then, by (\ref{eq:G}),
\[\begin{split}
G(z)&=1+\tr[m]n(u^5)+(wz)^\frac{1}{2}+\frac{\tr[m]n(u^5+u)z^4+\tr[m]n(u^5)w^2 z^2+\tr[m]n(u^4)z}{(z^2+wz+1)^2}\\
&\stackrel{(\ast)}{=}1+\tr[m]n(u^5)+(wz)^\frac{1}{2}+\frac{(w^5+w^3)z^4+w^3(1+w+w^2)^2z^2+w^4 z}{(z^2+wz+1)^2}\\
&=1+\tr[m]n(u^5)+(w^2+w^5+w^\frac{1}{2})f_0(z)\enspace,
\end{split}\]
where $(\ast)$ follows by $w(1+w+w^2)^2=\tr[m]n(u^5)$ and $f_0(z)$ is
an o-polynomial from (\ref{eq:m_arb}).

\section{Bent Functions from Adelaide Hyperovals}
Here we define o-polynomials that give rise to the Adelaide family of
hyperovals.

\begin{theorem}[Theorem~3.1~\cite{ChOKePe03}]
 \label{th:Adel}
Assume $m$ is even, $n=2m$ and denote $l=\frac{2^m-1}{3}$. Take any
$\beta\in\GF n\setminus\{1\}$ with $\beta^{2^m+1}=1$ and define the
following functions over $\GF m$
\[\begin{split}
f(x)&=\frac{\tr[m]n(\beta^l)(x+1)}{\tr[m]n(\beta)}
+\frac{\tr[m]n\big((\beta x+\beta^{-1})^l\big)}{\tr[m]n(\beta)(x+\tr[m]n(\beta)x^{1/2}+1)^{l-1}}
+x^{\frac{1}{2}}\quad\mbox{and}\\
eg(x)&=\frac{\tr[m]n(\beta^l)}{\tr[m]n(\beta)}x
+\frac{\tr[m]n\big((\beta^2 x+1)^l\big)}{\tr[m]n(\beta)\tr[m]n(\beta^l)(x+\tr[m]n(\beta)x^{1/2}+1)^{l-1}}
+\frac{1}{\tr[m]n(\beta^l)}x^{\frac{1}{2}}\enspace,
\end{split}\]
where
$e=\frac{\tr[m]n(\beta^l)}{\tr[m]n(\beta)}+\frac{1}{\tr[m]n(\beta^l)}+1$.
Then $g(x)$ and $f_s(x)$ (defined in (\ref{eq:fs})) are o-polynomials
for any $s\in\GF m$.
\end{theorem}

In particular, using that $\beta^{2^m}=\beta^{-1}$ we obtain that
\[e\tr[m]n(\beta)\tr[m]n(\beta^l)f_1(x)=\tr[m]n(\beta^{2l})
+\frac{\tr[m]n\big((x+\beta^2)^l\big)}{(x+\tr[m]n(\beta)x^{1/2}+1)^{l-1}}
+\tr[m]n(\beta)x^{\frac{1}{2}}\enspace.\]

For even $m$, take the following Niho bent function over $\GF n$
\[f(t)=\tr m(a t^{2^m+1})+\tr n(b t^{(2^m-1)\frac{1}{6}+1})\enspace,\]
where $\frac{1}{6}=\frac{2^{m-1}+1}{3}$ is an inverse of $6$ modulo
$2^m+1$, $a\in\GF m^*$ and $b\in\GF n^*$ are such that $b^{2^m+1}=a$.
As noted above, without loss of generality, it can be assumed that
$a=b=1$.

Assume $v=1$ and take $u\in\GF n\setminus\{1\}$ with $u^{2^m+1}=1$
that means $u\in\GF n\setminus\GF m$. Then $(u,1)$ is a basis of $\GF
n$ as a two-dimensional vector space over $\GF m$. Then for any
$x,y\in\GF m$, we obtain $f(ux+vy)$ having the form of
(\ref{eq:H_class}) with
\[\begin{split}
H(z)&=(z+u)^{\frac{2^m+1}{2}}+\tr[m]n\big((z+u)^{(2^m-1)\frac{1}{6}+1}\big)\\
\mu&=1\enspace.
\end{split}\]
Here all notation are from Subsection~\ref{subsec:class_H}.

Denote $d=(2^m-1)\frac{1}{6}+1=(2^{m-1}+1)l+1$, where
$l=\frac{2^m-1}{3}$. Then
\[2^{m+1}d\pmod{2^n-1}=(2^{m+1}+1)l+2^{m+1}=(2^m+1)(2l+1)+2l\]
and
\[\begin{split}
\tr[m]n\big((z+u)^{2d}\big)&=\tr[m]n\big((z+u)^{2^{m+1}d}\big)\\
&=(z+u)^{(2^m+1)(2l+1)}\tr[m]n\big((z+u)^{2l}\big)\\
&=(z^2+\tr[m]n(u)z+1)^{2l+1}\tr[m]n\big((z+u)^{2l}\big)\\
&=\frac{\tr[m]n\big((z+u)^{2l}\big)}{(z^2+\tr[m]n(u)z+1)^{l-1}}
\end{split}\]
since $3l=2^m-1$ and $z^2+\tr[m]n(u)z+1\in\GF m$.

Therefore, with $z\in\GF m$ and assuming $u=\beta^2$,
\[\begin{split}
G(z)&=1+\tr[m]n(\beta)z^{\frac{1}{2}}+\frac{\tr[m]n\big((z+\beta^2)^l\big)}{(z+\tr[m]n(\beta)z^{1/2}+1)^{l-1}}\\
&=1+\tr[m]n(\beta^{2l})+e\tr[m]n(\beta)\tr[m]n(\beta^l)f_1(z)\enspace.
\end{split}\]


\providecommand{\bysame}{\leavevmode\hbox
to3em{\hrulefill}\thinspace}
\providecommand{\MR}{\relax\ifhmode\unskip\space\fi MR }
\providecommand{\MRhref}[2]{%
  \href{http://www.ams.org/mathscinet-getitem?mr=#1}{#2}
} \providecommand{\href}[2]{#2}

\end{document}